\documentclass[thmsa,onecolumn,11pt,11pt]{article}
\usepackage{amssymb}
\usepackage{graphicx}
\usepackage{amsmath}

\newtheorem{theorem}{Theorem}

\newtheorem{corollary}[theorem]{Corollary}

\newtheorem{proposition}[theorem]{Proposition}

\newenvironment{proof}[1][Proof]{\textbf{#1.} }{\ \rule{0.5em}{0.5em}}

\begin{document}

\title{\textbf{Gaussian approximations of multiple integrals}}
\date{July 2, 2007}
\author{Giovanni PECCATI\thanks{%
Email: \textsf{giovanni.peccati@gmail.com}} \ \thanks{%
Partially supported by ISI Foundation, Lagrange Project.} \\
Laboratoire de Statistique Th\'{e}orique et Appliqu\'{e}e, Universit\'{e}
Paris VI}
\maketitle

\begin{abstract}
Fix $k\geq 1$, and let $\mathbf{I}\left( l\right) $, $l\geq 1$, be a
sequence of $k$-dimensional vectors of multiple Wiener-It\^{o} integrals
with respect to a general Gaussian process. We establish necessary and
sufficient conditions to have that, as $l\rightarrow +\infty $, the law of $%
\mathbf{I}\left( l\right) $ is asymptotically close (for example, in the
sense of Prokhorov's distance) to the law of a $k$-dimensional Gaussian
vector having the same covariance matrix as $\mathbf{I}\left( l\right) $.
The main feature of our results is that they require minimal assumptions
(basically, boundedness of variances) on the asymptotic behaviour of the
variances and covariances of the elements of $\mathbf{I}\left( l\right) $.
In particular, we will not assume that the covariance matrix of $\mathbf{I}%
\left( l\right) $ is convergent. This generalizes the results proved in
Nualart and Peccati (2005), Peccati and Tudor (2005) and Nualart and
Ortiz-Latorre (2007). As shown in Marinucci and Peccati (2007b), the
criteria established in this paper are crucial in the study of the
high-frequency behaviour of stationary fields defined on homogeneous spaces.


\end{abstract}

\section{Introduction}

Let $\mathbf{U}\left( l\right) =\left( U_{1}\left( l\right) ,...,U_{k}\left(
l\right) \right) $, $l\geq 1$, be a sequence of centered random observations
(not necessarily independent) with values in $\mathbb{R}^{k}$. Suppose that
the application $l \mapsto \mathbb{E}U_{i}\left( l\right) ^{2}$ is bounded
for every $i$, and also that the sequence of covariances $c_{l}\left(
i,j\right) =\mathbb{E}U_{i}\left( l\right) U_{j}\left( l\right) $ \textsl{%
does not }converge as $l\rightarrow +\infty $ (that is, for some fixed $%
i\neq j$, the limit $\lim_{l\rightarrow \infty }c_{l}\left( i,j\right) $
does not exist). Then, a natural question is the following: \textsl{is it
possible to establish criteria ensuring that, for large }$l$\textsl{, the
law of }$\mathbf{U}\left( l\right) $ \textsl{is close (in the sense of some
distance between probability measures) to the law of a Gaussian vector }$%
\mathbf{N}\left( l\right) =\left( N_{1}\left( l\right) ,...,N_{k}\left(
l\right) \right) $ \textsl{such that} $\mathbb{E}N_{i}\left( l\right)
N_{j}\left( l\right) =\mathbb{E}U_{i}\left( l\right) U_{j}\left( l\right)
=c_{l}\left( i,j\right) $\textsl{? }Note that the question is not trivial,
since the asymptotic irregularity of the covariance matrix $c_{l}\left(
\cdot ,\cdot \right) $ may in general prevent $\mathbf{U}\left( l\right) $
from converging in law toward a $k$-dimensional Gaussian distribution.

In this paper, we shall provide an exhaustive answer to the problem above in
the special case where the sequence $\mathbf{U}\left( l\right) $ has the
form
\begin{equation}
\mathbf{U}\left( l\right) =\mathbf{I}\left( l\right) =\left( I_{d_{1}}\left(
f_{l}^{\left( 1\right) }\right) ,...,I_{d_{k}}\left( f_{l}^{\left( k\right)
}\right) \right) \text{, \ \ }l\geq 1\text{,}  \label{l}
\end{equation}%
where the integers $d_{1},...,d_{k}\geq 1$ do not depend on $l$, $I_{d_{j}}$
indicates a multiple stochastic integral of order $d_{j}$ (with respect to
some isonormal Gaussian process $X$ over a Hilbert space $\mathfrak{H}$ --
see Section \ref{S : Pre} below for definitions), and each $f_{l}^{\left(
j\right) }\in \mathfrak{H}^{\odot d_{j}}$, $j=1,...,k$, is a symmetric
kernel. In particular, we shall prove that, whenever the elements of the
vectors $\mathbf{I}\left( l\right) $ have bounded variances (and without any
further requirements on the covariance matrix of $\mathbf{I}\left( l\right) $%
), the following three conditions are equivalent as $l\rightarrow +\infty $:

\begin{description}
\item[(\textbf{i})] $\gamma \left( \mathcal{L}\left( \mathbf{I}\left(
l\right) \right) ,\mathcal{L}\left( \mathbf{N}\left( l\right) \right)
\right) $ $\rightarrow 0$, where $\mathcal{L}\left( \cdot \right) $
indicates the law of a given random vector, $\mathbf{N}\left( l\right) $ is
a Gaussian vector having the same covariance matrix as $\mathbf{I}\left(
l\right) $, and $\gamma $ is some appropriate metric on the space of
probability measures on $\mathbb{R}^{k}$;

\item[(\textbf{ii})] For every $j=1,...,k,$ $\mathbb{E}\left(
I_{d_{j}}\left( f_{l}^{\left( j\right) }\right) ^{4}\right) -3\mathbb{E}%
\left( I_{d_{j}}\left( f_{l}^{\left( j\right) }\right) ^{2}\right)
^{2}\rightarrow 0;$

\item[(\textbf{iii})] For every $j=1,...,k$ and every $p=1,...,d_{j}-1$, the
sequence of contractions (to be formally defined in Section \ref{S : Pre}) $%
f_{l}^{\left( j\right) }\otimes _{p}f_{l}^{\left( j\right) }$, $l\geq 1$, is
such that
\begin{equation}
f_{l}^{\left( j\right) }\otimes _{p}f_{l}^{\left( j\right) }\rightarrow 0%
\text{ \ \ in \ \ }\mathfrak{H}^{\odot 2\left( d_{j}-p\right) }\text{.}
\label{qq}
\end{equation}
\end{description}

Some other conditions, involving for instance Malliavin operators, are
derived in the subsequent sections. As discussed in Section \ref{S : Concl},
our results are motivated by the derivation of high-frequency Gaussian
approximations of stationary fields defined on homogeneous spaces -- a
problem tackled in \cite{MaPe} and \cite{MaPe2}.

Note that the results of this paper are a generalization of the following
theorem, which combines results proved in \cite{NuOL}, \cite{NuPe} and \cite%
{PT}.

\bigskip

\textbf{Theorem 0. }\textsl{Suppose that the vector }$\mathbf{I}\left(
l\right) $\textsl{\ in (\ref{l}) is such that, as }$l\rightarrow +\infty $%
\textsl{,}%
\begin{equation*}
\mathbb{E}I_{d_{i}}\left( f_{l}^{\left( i\right) }\right) I_{d_{j}}\left(
f_{l}^{\left( j\right) }\right) \rightarrow \mathbf{C}\left( i,j\right)
\text{, \ \ }1\leq i,j\leq k\text{,}
\end{equation*}%
\textsl{where }$\mathbf{C}=\left\{ \mathbf{C}\left( i,j\right) \right\} $%
\textsl{\ is some positive definite matrix. Then, the following four
conditions are equivalent, as }$l\rightarrow +\infty $\textsl{:}

\begin{enumerate}
\item $\mathbf{I}\left( l\right) \overset{Law}{\rightarrow }\mathbf{N}\left(
0,\mathbf{C}\right) $, \textsl{where }$\mathbf{N}\left( 0,\mathbf{C}\right) $%
\textsl{\ is a }$k$\textsl{-dimensional centered Gaussian vector with
covariance matrix }$\mathbf{C}$\textsl{;}

\item \textsl{Relation (\ref{qq}) takes place for every }$j=1,...,k$\textsl{%
\ and every }$p=1,...,d_{j}-1$\textsl{;}

\item \textsl{For every }$j=1,...k$\textsl{, }$\mathbb{E}\left(
I_{d_{j}}\left( f_{l}^{\left( j\right) }\right) ^{4}\right)\rightarrow3%
\mathbf{C}\left( j,j\right) ^{2}$\textsl{;}

\item \textsl{For every }$j=1,...k$\textsl{, }$\left\Vert D\left[
I_{d_{j}}\left( f_{l}^{\left( j\right) }\right) \right] \right\Vert _{%
\mathfrak{H}}^{2}\rightarrow d_{j} $\textsl{\ \ in }$L^{2}$\textsl{, where }$%
D\left[ I_{d_{j}}\left( f_{l}^{\left( j\right) }\right) \right] $\textsl{\
denotes the Malliavin derivative of }$I_{d_{j}}\left( f_{l}^{\left( j\right)
}\right) $\textsl{\ (see the next section).}
\end{enumerate}

\bigskip

The equivalence of Points 1.-3. in the case $k=1$ has been first proved in
\cite{NuPe} by means of the Dambis-Dubins-Schwarz (DDS) Theorem (see \cite[%
Ch.\ V]{RY}), whereas the proof in the case $k\geq 2$ has been achieved (by
similar techniques) in \cite{PT}; the fact that Point 4. is also necessary
and sufficient for the CLT\ at Point 1. has been recently proved in \cite%
{NuOL}, by means of a Malliavin calculus approach. For some applications of
Theorem 0 (in quite different frameworks), see e.g. \cite{CorNuWoe}, \cite%
{DPY}, \cite{HuNu}, \cite{MaPe} or \cite{NeuNourdin}.

The techniques we use to achieve our main results are once again the DDS
Theorem, combined with Burkholder-Davis-Gundy inequalities and some results
(taken from \cite[Section 11.7]{Dudley book}) concerning `uniformities' over
classes of probability measures.

The paper is organized as follows. In Section \ref{S : Pre} we discuss some
preliminary notions concerning Gaussian fields, multiple integrals and
metrics on probabilities. Section \ref{S :Main} contains the statements of
the main results of the paper. The proof of Theorem \ref{T : Appendix} (one
of the crucial results of this note) is achieved in Section \ref{S : Proof}.
Section \ref{S : Concl} is devoted to applications.

\section{Preliminairies\label{S : Pre}}

We present a brief review of the main notions and results that are needed in
the subsequent sections. The reader is referred to \cite{Janss} or \cite[Ch.
1]{NualartBook} for any unexplained definition.

\begin{description}
\item[Hilbert spaces.] In what follows, the symbol $\mathfrak{H}$ indicates
a real separable Hilbert space, with inner product $\left\langle \cdot
,\cdot \right\rangle _{\mathfrak{H}}$ and norm $\left\Vert \cdot \right\Vert
_{\mathfrak{H}}$. For every $d\geq 2$, we denote by $\mathfrak{H}^{\otimes
2} $ and $\mathfrak{H}^{\odot 2}$, respectively, the $n$th \textsl{tensor
product} of $\mathfrak{H}$ and the $n$th \textsl{symmetric tensor product}
of $\mathfrak{H}$. We also write $\mathfrak{H}^{\otimes 1}$ $=\mathfrak{H}%
^{\odot 1}$ $=\mathfrak{H}$.

\item[Isonormal Gaussian processes.] We write $X=\left\{ X\left( h\right)
:h\in \mathfrak{H}\right\} $ to indicate an \textsl{isonormal Gaussian
process} over $\mathfrak{H}$. This means that $X$ is a collection of
real-valued, centered and (jointly) Gaussian random variables indexed by the
elements of $\mathfrak{H}$, defined on a probability space $\left( \Omega ,%
\mathcal{F},\mathbb{P}\right) $ and such that, for every $h,h^{\prime }\in
\mathfrak{H}$,
\begin{equation*}
\mathbb{E}\left[ X\left( h\right) X\left( h^{\prime }\right) \right]
=\left\langle h,h^{\prime }\right\rangle _{\mathfrak{H}}\text{.}
\end{equation*}%
We denote by $L^{2}\left( X\right) $ the (Hilbert) space of the real-valued
and square-integrable functionals of $X$.

\item[Isometry, chaoses and multiple integrals.] For every $d\geq 1$ we will
denote by $I_{d}$ the isometry between $\mathfrak{H}^{\odot d}$ equipped
with the norm $\sqrt{d!}\left\Vert \cdot \right\Vert _{\mathfrak{H}^{\otimes
d}}$ and the $d$th Wiener chaos of $X$. In the particular case where $%
\mathfrak{H}=L^{2}\left( A,\mathcal{A},\mu \right) $, $\left( A,\mathcal{A}%
\right) $ is a measurable space, and $\mu $ is a $\sigma $-finite and
non-atomic measure, then $\mathfrak{H}^{\odot d}=L_{s}^{2}\left( A^{d},%
\mathcal{A}^{\otimes d},\mu ^{\otimes d}\right) $ is the space of symmetric
and square integrable functions on $A^{d}$ and for every $f\in \mathfrak{H}%
^{\odot d}$, $I_{d}\left( f\right) $ is the \textsl{multiple Wiener-It\^{o}
integral} (of order $d$) of $f$ with respect to $X$, as defined e.g. in \cite%
[Ch. 1]{NualartBook}. It is well-known that a random variable of the type $%
I_d(f)$, where $d\geq2$ and $f\neq0$, \textsl{cannot be Gaussian}. Moreover,
every $F\in L^{2}\left( X\right) $ admits a unique \textsl{Wiener}\ \textsl{%
chaotic decomposition} of the type $F=\mathbb{E}\left( F\right) $ $%
+\sum_{d=1}^{\infty }I_{d}\left( f_{d}\right) $, where $f_{d}\in \mathfrak{H}%
^{\odot d}$, $d\geq 1$, and the convergence of the series is in $L^{2}\left(
X\right) $.

\item[Malliavin derivatives.] We will use Malliavin derivatives in Section %
\ref{S :Main}, where we generalize some of the results proved in \cite{NuOL}%
. The class $\mathcal{S}$ of \textsl{smooth} random variables is defined as
the collection of all functionals of the type%
\begin{equation}
F=f\left( X\left( h_{1}\right) ,...,X\left( h_{m}\right) \right) ,  \label{D}
\end{equation}%
where $h_{1},...,h_{m}\in \mathfrak{H}$ and $f$ is bounded and has bounded
derivatives of all order. The operator $D$, called the \textsl{Malliavin
derivative operator}, is defined on $\mathcal{S}$ by the relation%
\begin{equation*}
DF=\sum_{i=1}^{M}\frac{\partial }{\partial x_{i}}f\left(
h_{1},...,h_{m}\right) h_{i}\text{,}
\end{equation*}%
where $F$ has the form (\ref{D}). Note that $DF$ is an element of $%
L^{2}\left( \Omega ;\mathfrak{H}\right) $. As usual, we define the domain of
$D$, noted $\mathbb{D}^{1,2}$, to be the closure of $\mathcal{S}$ with
respect to the norm $\left\Vert F\right\Vert _{1,2}\triangleq \mathbb{E}%
\left( F^{2}\right) $ $+\mathbb{E}\left\Vert DF\right\Vert _{\mathfrak{H}%
}^{2}$. When $F\in \mathbb{D}^{1,2}$, we may sometimes write $DF=D\left[ F%
\right] $, depending on the notational convenience. Note that any finite sum
of multiple Wiener-It\^{o} integrals is an element of $\mathbb{D}^{1,2}$.

\item[Contractions.] Let $\{e_{k}:k\geq 1\}$ be a complete orthonormal
system of $\mathfrak{H}$. For any fixed $f\in \mathfrak{H}^{\odot n}$, $g\in
\mathfrak{H}^{\odot m}$ and $p\in \left\{ 0,...,n\wedge m\right\} $, we
define the $p$th \textsl{contraction}\textit{\ }of $f$ and $g$ to be the
element of $\mathfrak{H}^{\otimes n+m-2p}$ given by
\begin{equation*}
f\otimes _{p}g=\sum_{i_{1},\ldots ,i_{p}=1}^{\infty }\ \left\langle
f,e_{i_{1}}\otimes \cdots \otimes e_{i_{p}}\right\rangle _{\mathfrak{H}%
^{\otimes p}}\otimes \left\langle g,e_{i_{1}}\otimes \cdots \otimes
e_{i_{p}}\right\rangle _{\mathfrak{H}^{\otimes p}}\text{.}
\end{equation*}%
We stress that $f\otimes _{p}g$ need not be an element of $\mathfrak{H}%
^{\odot n+m-2p}$.\ We denote by $f\widetilde{\otimes }_{p}g$ the
symmetrization of $f\otimes _{p}g$. Note that $f\otimes _{0}g$ is just the
tensor product $f\otimes g$ of $f$ and $g$. If $n=m$, then $f\otimes
_{n}g=\left\langle f,g\right\rangle _{\mathfrak{H}^{\otimes n}}$.

\item[Metrics on probabilities.] For $k\geq 1$ we define $\mathbf{P}\left(
\mathbb{R}^{k}\right) $ to be the class of all probability measures on $%
\mathbb{R}^{k}$. Given a metric $\gamma \left( \cdot ,\cdot \right) $ on $%
\mathbf{P}\left( \mathbb{R}^{k}\right) $, we say that $\gamma $ \textsl{%
metrizes the weak convergence} \textsl{on} $\mathbf{P}\left( \mathbb{R}%
^{k}\right) $ whenever the following double implication holds for every $%
Q\in \mathbf{P}\left( \mathbb{R}^{k}\right) $ and every $\left\{ Q_{l}:l\geq
1\right\} \subset \mathbf{P}\left( \mathbb{R}^{k}\right) $ (as $l\rightarrow
+\infty $): $\gamma \left( Q_{l},Q\right) \rightarrow 0$ if, and only if, $%
Q_{l}$ converges weakly to $Q$. Some examples of metrizing $\gamma $ are the
\textsl{Prokhorov metric }(usually noted $\rho $) or the \textsl{%
Fortet-Mounier metric }(usually noted $\beta $).\textsl{\ }Recall that
\begin{equation}
\rho \left( P,Q\right) =\inf \{\epsilon >0:P\left( A\right) \leq
Q(A^{\epsilon })+\epsilon ,\text{ \ for every Borel set }A\subset \mathbb{R}%
^{k}\}  \label{PR met}
\end{equation}%
where $A^{\epsilon }=\{x:\left\Vert x-y\right\Vert <\varepsilon $ \ for some
$y\in A\}$, and $\left\Vert \cdot \right\Vert $ is the Euclidiean norm. Also,%
\begin{equation}
\beta \left( P,Q\right) = \sup \left\{ \left\vert \int fd\left(
P-Q\right) \right\vert :\left\Vert f\right\Vert _{BL}\leq
1\right\} \text{,} \label{FM met}
\end{equation}%
where $\left\Vert \cdot \right\Vert _{BL}=\left\Vert \cdot \right\Vert _{L}$
$+\left\Vert \cdot \right\Vert _{\infty }$, and $\left\Vert \cdot
\right\Vert _{L}$ is the usual Lipschitz seminorm (see \cite[p. 394]{Dudley
book} for further details). The fact that we focus on the Prokhorov and the
Fortet-Mounier metric is due to the following fact, proved in \cite[Th.\
11.7.1]{Dudley book}. \textsl{For any two sequences }$\left\{ P_{l}\right\}
,\left\{ Q_{l}\right\} \subset \mathbf{P}\left( \mathbb{R}^{k}\right) $%
\textsl{, the following three conditions }(\textbf{A})--(\textbf{C})\textsl{%
\ are equivalent: }\textbf{(A)}\textsl{\ }$\lim_{l\rightarrow +\infty }\beta
\left( P_{l},Q_{l}\right) $\textsl{\ }$=$\textsl{\ }$0$\textsl{; }\textbf{(B)%
}\textsl{\ }$\lim_{l\rightarrow +\infty }\rho \left( P_{l},Q_{l}\right) =0$%
\textsl{; }\textbf{(C)}\textsl{\ on some auxiliary probability space }$%
\left( \Omega ^{\ast },\mathcal{F}^{\ast },\mathbb{P}^{\ast }\right) $%
\textsl{, there exist sequences of random vectors }$\left\{ \mathbf{N}^{\ast
}\left( l\right) :l\geq 1\right\} $\textsl{\ and }$\left\{ \mathbf{I}^{\ast
}\left( l\right) :l\geq 1\right\} $\textsl{\ such that }%
\begin{equation}
\mathcal{L}\left( \mathbf{I}^{\ast }\left( l\right) \right) =P_{l}\text{
\textsl{and}\ }\mathcal{L}\left( \mathbf{N}^{\ast }\left( l\right) \right)
=Q_{l}\text{ \ \textsl{for every} }l\text{, \textsl{and} }\left\Vert \mathbf{%
I}^{\ast }\left( l\right) -\mathbf{N}^{\ast }\left( l\right) \right\Vert
\rightarrow 0,\text{ \ \textsl{a.s.}-}\mathbb{P}^{\ast }\text{,}
\label{feff}
\end{equation}%
\textsl{where }$\mathcal{L}\left( \mathbf{\cdot }\right) $\textsl{\
indicates the law of a given random vector, and }$\left\Vert \cdot
\right\Vert $\textsl{\ is the Euclidean norm.}
\end{description}

\section{Main results\label{S :Main}}

Fix integers $k\geq 1$ and $d_{1},...,d_{k}\geq 1$, and consider a sequence
of $k$-dimensional random vectors of the type%
\begin{equation}
\mathbf{I}\left( l\right) =\left( I_{d_{1}}\left( f_{l}^{\left( 1\right)
}\right) ,...,I_{d_{k}}\left( f_{l}^{\left( k\right) }\right) \right) \text{
, \ \ }l\geq 1\text{ ,}  \label{Al-sequence}
\end{equation}%
where, for each $l\geq 1$ and every $j=1,...,k$, $f_{l}^{\left( j\right) }$
is an element of $\mathfrak{H}^{\odot d_{j}}$. We will suppose the following:

\begin{itemize}
\item There exists $\eta >0$ such that $\left\Vert f_{l}^{\left( j\right)
}\right\Vert _{\mathfrak{H}^{\odot d_{j}}}\geq \eta $, for every $j=1,...,k$
and every $l\geq 1$.

\item For every $j=1,...,k$, the sequence
\begin{equation}
\mathbb{E}\left[ I_{d_{j}}\left( f_{l}^{\left( j\right) }\right) ^{2}\right]
=d_{j}!\left\Vert f_{l}^{\left( j\right) }\right\Vert _{\mathfrak{H}^{\odot
d_{j}}}^{2}\text{, \ \ }l\geq 1\text{,}  \label{Avarbound}
\end{equation}%
is bounded.
\end{itemize}

Note that the integers $d_{1},...,d_{k}$ do not depend on $l$. For every $%
l\geq 1$, we denote by $\mathbf{N}\left( l\right) =\left( N_{l}^{\left(
1\right) },...,N_{l}^{\left( k\right) }\right) $ a centered $k$-dimensional
Gaussian vector with the same covariance matrix as $\mathbf{I}\left(
l\right) $, that is,%
\begin{equation}
\mathbb{E}\left[ N_{l}^{\left( i\right) }N_{l}^{\left( j\right) }\right] =%
\mathbb{E}\left[ I_{d_{i}}\left( f_{l}^{\left( i\right) }\right)
I_{d_{j}}\left( f_{l}^{\left( j\right) }\right) \right] \text{,}
\label{AGcov}
\end{equation}%
for every $1\leq i,j\leq k$. For every $\mathbf{\lambda }=\left( \lambda
_{1},...,\lambda _{k}\right) \in \mathbb{R}^{k}$, we also use the compact
notation: $\left\langle \mathbf{\lambda },\mathbf{I}\left( l\right)
\right\rangle _{k}=$ $\sum_{j=1}^{k}\lambda _{j}I_{d_{j}}\left(
f_{l}^{\left( j\right) }\right) $ and $\left\langle \mathbf{\lambda },%
\mathbf{N}\left( l\right) \right\rangle _{k}=$ $\sum_{j=1}^{k}\lambda
_{j}N_{l}^{\left( j\right) }$.

The next result is one of the main contributions of this paper. Its proof is
deferred to Section \ref{S : Proof}.

\begin{theorem}
\label{T : Appendix}Let the above notation and assumptions prevail, and
suppose that, for every $j=1,...,k$, the following asymptotic condition
holds: for every $p=1,...,d_{j}-1$,
\begin{equation}
\left\Vert f_{l}^{\left( j\right) }\otimes _{p}f_{l}^{\left( j\right)
}\right\Vert _{\mathfrak{H}^{\odot 2\left( d_{j}-p\right) }}\rightarrow 0%
\text{, \ \ as }l\rightarrow +\infty \text{.}  \label{Acontrconv}
\end{equation}%
Then, as $l\rightarrow +\infty $ and for every compact set $M\subset \mathbb{%
R}^{k}$,
\begin{equation}
\sup_{\mathbf{\lambda }\in M}\left\vert \mathbb{E}\left[ \exp \left( \mathrm{%
i}\left\langle \mathbf{\lambda },\mathbf{I}\left( l\right) \right\rangle
_{k}\right) \right] -\mathbb{E}\left[ \exp \left( \mathrm{i}\left\langle
\mathbf{\lambda },\mathbf{N}\left( l\right) \right\rangle _{k}\right) \right]
\right\vert \rightarrow 0\text{.}  \label{AunifFour}
\end{equation}
\end{theorem}

We now state two crucial consequences of Theorem \ref{T :
Appendix}. The first one (Proposition \ref{C : Appendix}) provides
a formal meaning to the intuitive fact that, since
(\ref{AunifFour}) holds and since the
variances of $\mathbf{I}\left( l\right) $ do not explode, the laws of $%
\mathbf{I}\left( l\right) $ and $\mathbf{N}\left( l\right) $ are
\textquotedblleft asymptotically close\textquotedblright . The second one
(Theorem \ref{T : Gen}) combines Theorem \ref{T : Appendix} and Proposition %
\ref{C : Appendix} to obtain an exhaustive generalization \textquotedblleft
without covariance conditions\textquotedblright\ of Theorem 0 (see the
Introduction). Note that in the statement of Theorem \ref{T : Gen} also
appear Malliavin operators, so that our results are a genuine extension of
the main findings by Nualart and Ortiz-Latorre in \cite{NuOL}. We stress
that multiple stochastic integrals of the type $I_{d}\left( f\right) $, $%
d\geq 1$ and $f\in \mathfrak{H}^{\odot d}$, are always such that $%
I_{d}\left( f\right) \in \mathbb{D}^{1,2}.$

\begin{proposition}
\label{C : Appendix}Let the assumptions of Theorem \ref{T : Appendix}
prevail (in particular, (\ref{Acontrconv}) holds), and denote by $\mathcal{L}%
\left( \mathbf{I}\left( l\right) \right) $ and $\mathcal{L}\left( \mathbf{N}%
\left( l\right) \right) $, respectively, the law of $\mathbf{I}\left(
l\right) $ and $\mathbf{N}\left( l\right) $, $l\geq 1$. Then, the two
collections $\left\{ \mathcal{L}\left( \mathbf{N}\left( l\right) \right)
:l\geq 1\right\} $ and $\left\{ \mathcal{L}\left( \mathbf{I}\left( l\right)
\right) :l\geq 1\right\} $ are tight. Moreover, if $\gamma \left( \cdot
,\cdot \right) $ metrizes the weak convergence on $\mathbf{P}\left( \mathbb{R%
}^{k}\right) $, then
\begin{equation}
\lim_{l\rightarrow +\infty }\gamma \left( \mathcal{L}\left( \mathbf{I}\left(
l\right) \right) ,\mathcal{L}\left( \mathbf{N}\left( l\right) \right)
\right) =0\text{.}  \label{Aconvgamma}
\end{equation}
\end{proposition}

\begin{proof}
The fact that $\left\{ \mathcal{L}\left( \mathbf{N}\left( l\right) \right)
:l\geq 1\right\} $ and $\left\{ \mathcal{L}\left( \mathbf{I}\left( l\right)
\right) :l\geq 1\right\} $ are tight is a consequence of the boundedness of
the sequence (\ref{Avarbound}) and of the relation $\mathbb{E[}%
I_{d_{j}}(f_{l}^{\left( j\right) })^{2}]=\mathbb{E[}(N_{l}^{\left( j\right)
})^{2}]$. The rest of the proof is standard, and is provided for the sake of
completeness. We shall prove (\ref{Aconvgamma}) by contradiction. Suppose
there exist $\varepsilon >0$ and a subsequence $\left\{ l_{n}\right\} $ such
that $\gamma \left( \mathcal{L}\left( \mathbf{I}\left( l_{n}\right) \right) ,%
\mathcal{L}\left( \mathbf{N}\left( l_{n}\right) \right) \right) >\varepsilon
$ for every $n$. Tightness implies that $\left\{ l_{n}\right\} $ must
contain a subsequence $\left\{ l_{n^{\prime }}\right\} $ such that $\mathcal{%
L}\left( \mathbf{I}\left( l_{n^{\prime }}\right) \right) $ and $\mathcal{L}%
\left( \mathbf{N}\left( l_{n^{\prime }}\right) \right) $ are both weakly
convergent. Since (\ref{AunifFour}) holds, we deduce that $\mathcal{L}\left(
\mathbf{I}\left( l_{n^{\prime }}\right) \right) $ and $\mathcal{L}\left(
\mathbf{N}\left( l_{n^{\prime }}\right) \right) $ must necessarily converge
to the same weak limit, say $Q$. The fact that $\gamma $ metrizes the weak
convergence implies finally that
\begin{equation}
\gamma \left( \mathcal{L}\left( \mathbf{I}\left( l_{n^{\prime }}\right)
\right) ,\mathcal{L}\left( \mathbf{N}\left( l_{n^{\prime }}\right) \right)
\right) \leq \gamma \left( \mathcal{L}\left( \mathbf{I}\left( l_{n^{\prime
}}\right) \right) ,Q\right) +\gamma \left( \mathcal{L}\left( \mathbf{N}%
\left( l_{n^{\prime }}\right) \right) ,Q\right) \underset{n^{\prime
}\rightarrow +\infty }{\rightarrow }0\text{,}  \label{Ain}
\end{equation}%
thus contradicting the former assumptions on $\left\{ l_{n}\right\} $ (note
that the inequality in (\ref{Ain}) is just the triangle inequality). This
shows that (\ref{Aconvgamma}) must necessarily take place.
\end{proof}

\bigskip

\textbf{Remarks. }(i) A result analogous to the arguments used in the proof
of Corollary \ref{C : Appendix} is stated in \cite[Exercise 3, p. 419]%
{Dudley book}. Note also that, without tightness, a condition such as (\ref%
{AunifFour}) \textsl{does not allow }to deduce the asymptotic relation (\ref%
{Aconvgamma}). See for instance \cite[Proposition 11.7.6]{Dudley book} for a
counterexample involving the Prokhorov metric on $\mathbf{P(}\mathbb{R)}$.

(ii) Since (\ref{Aconvgamma}) holds in particular when $\gamma $ is equal to
the Prokhorov metric or the Fortet-Mounier metric (as defined in (\ref{PR
met}) and (\ref{FM met})), Proposition \ref{C : Appendix} implies that, on
some auxiliary probability space $\left( \Omega ^{\ast },\mathcal{F}^{\ast },%
\mathbb{P}^{\ast }\right) $, there exist sequences of random vectors $%
\left\{ \mathbf{N}^{\ast }\left( l\right) :l\geq 1\right\} $ and $\left\{
\mathbf{I}^{\ast }\left( l\right) :l\geq 1\right\} $ such that
\begin{equation}
\mathbf{I}^{\ast }\left( l\right) \overset{law}{=}\mathbf{I}\left( l\right)
\text{ and\ }\mathbf{N}^{\ast }\left( l\right) \overset{law}{=}\mathbf{N}%
\left( l\right) \text{ for every }l\text{, and }\left\Vert \mathbf{I}^{\ast
}\left( l\right) -\mathbf{N}^{\ast }\left( l\right) \right\Vert \rightarrow
0,\text{ \ a.s.-}\mathbb{P}^{\ast }\text{,}  \label{reff}
\end{equation}%
where $\left\Vert \cdot \right\Vert $ stands for the Euclidean norm (see (%
\ref{feff}), as well as \cite[Theorem 11.7.1]{Dudley book}).

\begin{theorem}
\label{T : Gen}Suppose that the sequence $\mathbf{I}\left( l\right) $, $%
l\geq 1$, verifies the assumptions of this section (in particular, for every
$j=1,...,k$, the sequence of variances appearing in (\ref{Avarbound}) is
bounded). Then, the following conditions are equivalent.

\begin{enumerate}
\item As $l\rightarrow +\infty $, relation (\ref{Acontrconv}) is satisfied
for every $j=1,...,k$ and every $p=1,...,d_{j}-1$;

\item
\begin{equation}
\lim_{l\rightarrow +\infty }\rho \left( \mathcal{L}\left( \mathbf{I}\left(
l\right) \right) ,\mathcal{L}\left( \mathbf{N}\left( l\right) \right)
\right) =\lim_{l\rightarrow +\infty }\beta \left( \mathcal{L}\left( \mathbf{I%
}\left( l\right) \right) ,\mathcal{L}\left( \mathbf{N}\left( l\right)
\right) \right) =0  \label{bb}
\end{equation}%
where $\rho $ and $\beta $ are, respectively, the Prokhorov metric and the
Fortet-Mounier metric, as defined in (\ref{PR met}) and (\ref{FM met});

\item As $l\rightarrow +\infty $, for every $j=1,...,k$,
\begin{equation*}
\mathbb{E}\left[ I_{d_{j}}\left( f_{l}^{\left( j\right) }\right) ^{4}\right]
-3\mathbb{E}\left[ I_{d_{j}}\left( f_{l}^{\left( j\right) }\right) ^{2}%
\right] ^{2}=\mathbb{E}\left[ I_{d_{j}}\left( f_{l}^{\left(
j\right) }\right) ^{4}\right] -3(d_{j}!)^{2}\left\Vert
f_{l}^{\left( j\right) }\right\Vert _{\mathfrak{H}^{\odot
d_{j}}}^{4}\rightarrow 0;
\end{equation*}

\item For every $j=1,...,k$,
\begin{equation}
\lim_{l\rightarrow +\infty }\rho \left( \mathcal{L}\left( I_{d_{j}}\left(
f_{l}^{\left( j\right) }\right) \right) ,\mathcal{L}\left( N_{l}^{\left(
j\right) }\right) \right) =\lim_{l\rightarrow +\infty }\beta \left( \mathcal{%
L}\left( I_{d_{j}}\left( f_{l}^{\left( j\right) }\right) \right) ,\mathcal{L}%
\left( N_{l}^{\left( j\right) }\right) \right) =0,  \label{s-conv}
\end{equation}%
where $\rho $ and $\beta $ are the Prokhorov and Fortet-Mounier metric on $%
\mathbb{R}$;

\item For every $j=1,...,k$,%
\begin{equation}
\left\Vert D\left[ I_{d_{j}}\left( f_{l}^{\left( j\right) }\right) \right]
\right\Vert _{\mathfrak{H}}^{2}-d_{j}\left( d_{j}!\right) \left\Vert
f_{l}^{\left( j\right) }\right\Vert _{\mathfrak{H}^{\otimes
d_{j}}}^{2}\rightarrow 0\text{, \ \ in }L^{2}\left( X\right) \text{,}
\label{Mall}
\end{equation}%
as $l\rightarrow +\infty $, where $D$ is the Malliavin derivative operator
defined in Section \ref{S : Pre}.
\end{enumerate}
\end{theorem}

\begin{proof}
The implication 1. $\Longrightarrow $ 2., is a consequence of Theorem \ref{T
: Appendix} and Proposition \ref{C : Appendix}. Now suppose (\ref{bb}) is in
order. Then, according to \cite[Theorem 11.7.1]{Dudley book}, on a
probability space $\left( \Omega ^{\ast },\mathcal{F}^{\ast },\mathbb{P}%
^{\ast }\right) $, there exist sequences of random vectors $\mathbf{N}^{\ast
}\left( l\right) =(N_{l}^{\ast ,\left( 1\right) },...,N_{l}^{\ast ,\left(
j\right) })$, $l\geq 1$, and $\mathbf{I}^{\ast }\left( l\right) =\left(
I_{l}^{\ast ,\left( 1\right) },...,I_{l}^{\ast ,\left( k\right) }\right) $, $%
l\geq 1$, such that (\ref{reff}) takes place. Now%
\begin{equation*}
3\mathbb{E}\left[ I_{d_{j}}\left( f_{l}^{\left( j\right) }\right) ^{2}\right]
^{2}=3\mathbb{E}\left[ \left( N_{l}^{\left( j\right) }\right) ^{2}\right]
^{2}=\mathbb{E}\left[ \left( N_{l}^{\left( j\right) }\right) ^{4}\right] =%
\mathbb{E}^{\ast }\left[ \left( N_{l}^{\ast ,\left( j\right) }\right) ^{4}%
\right] \text{,}
\end{equation*}%
for every $j=1,...,k$, so that
\begin{equation}
\mathbb{E}\left[ I_{d_{j}}\left( f_{l}^{\left( j\right) }\right) ^{4}\right]
-3\mathbb{E}\left[ I_{d_{j}}\left( f_{l}^{\left( j\right) }\right) ^{2}%
\right] ^{2}=\mathbb{E}^{\ast }\left[ \left( I_{l}^{\ast ,\left( j\right)
}\right) ^{4}-\left( N_{l}^{\ast ,\left( j\right) }\right) ^{4}\right]
\underset{l\rightarrow +\infty }{\rightarrow }0\text{.}  \label{H}
\end{equation}%
The convergence to zero in (\ref{H}) is a consequence of the
boundedness of the sequence (\ref{Avarbound}), implying that the
family $A_{l}^{\ast }=\left( I_{l}^{\ast ,\left( j\right) }\right)
^{4}-\left( N_{l}^{\ast ,\left( j\right) }\right) ^{4}$, $l\geq
1$, is uniformly integrable. To see why $\{A_{l}^{\ast }\}$ is
uniformly integrable, one can use the fact that, since each
$I_{l}^{\ast ,\left( j\right) }$
has the same law as an element of the $d_{j}$th chaos of $X$ and each $%
N_{l}^{\ast ,\left( j\right) }$ is Gaussian, then (see e.g.
\cite[Ch. VI]{Janss}) for every $p\geq 2$ there
exists a universal positive constant $C_{p,j}$ (independent of $l$) such that%
\begin{eqnarray*}
\mathbb{E}\left[ \left\vert A_{l}^{\ast }\right\vert ^{p}\right] ^{1/p} &=&%
\mathbb{E}^{\ast }\left[ \left\vert \left( I_{l}^{\ast ,\left( j\right)
}\right) ^{4}-\left( N_{l}^{\ast ,\left( j\right) }\right) ^{4}\right\vert
^{p}\right] ^{1/p} \\
&\leq &\mathbb{E}^{\ast }\left[ \left( I_{l}^{\ast ,\left( j\right) }\right)
^{4p}\right] ^{4/4p}+\mathbb{E}\left[ \left( N_{l}^{\ast ,\left( j\right)
}\right) ^{4p}\right] ^{4/4p} \\
&\leq &C_{p,j}\mathbb{E}^{\ast }\left[ \left( I_{l}^{\ast ,\left( j\right)
}\right) ^{2}\right] ^{2}+C_{p,j}\mathbb{E}^{\ast }\left[ \left( N_{l}^{\ast
,\left( j\right) }\right) ^{2}\right] ^{2} \\
&=&2C_{p,j}\times \left( d_{j}!\right) ^{2}\left\Vert f_{l}^{\left( j\right)
}\right\Vert _{\mathfrak{H}^{\odot d_{j}}}^{4}\leq 2C_{p,j}M_{j}\text{,}
\end{eqnarray*}%
where $M_{j}=\sup_{l}\left( d_{j}!\right) ^{2}\left\Vert f_{l}^{\left(
j\right) }\right\Vert ^{4}<+\infty $, due to (\ref{Avarbound}). This proves
that 2. $\Longrightarrow $ 3.. The implication 3. $\Longrightarrow $ 1. can
be deduced from the formula (proved in \cite[p. 183]{NuPe})
\begin{eqnarray*}
&&\mathbb{E}\left[ I_{d_{j}}\left( f_{l}^{\left( j\right) }\right) ^{4}%
\right] -3\mathbb{E}\left[ I_{d_{j}}\left( f_{l}^{\left( j\right) }\right)
^{2}\right] ^{2}=\mathbb{E}\left[ I_{d_{j}}\left( f_{l}^{\left( j\right)
}\right) ^{4}\right] -3\left( d_{j}!\right) ^{2}\left\Vert f_{l}^{\left(
j\right) }\right\Vert _{\mathfrak{H}^{\otimes d_{j}}}^{4} \\
&=&\sum_{p=1}^{d_{j}-1}\frac{\left( d_{j}!\right) ^{4}}{\left( p!\left(
d_{j}-p\right) !\right) ^{2}}\left\{ \left\Vert f_{l}^{\left( j\right)
}\otimes _{p}f_{l}^{\left( j\right) }\right\Vert _{\mathfrak{H}^{\otimes
2\left( d_{j}-p\right) }}^{2}+\binom{2\left( d_{j}-p\right) }{d_{j}-p}%
\left\Vert f_{l}^{\left( j\right) }\widetilde{\otimes }_{p}f_{l}^{\left(
j\right) }\right\Vert _{\mathfrak{H}^{\otimes 2\left( d_{j}-p\right)
}}^{2}\right\} ,
\end{eqnarray*}

The equivalence 1. $\Longleftrightarrow $ 4. is an immediate consequence of
the previous discussion.

To conclude the proof, we shall now show the double implication 1. $%
\Longleftrightarrow $ 5.. To do this, we first observe that, by performing
the same caclulations as in \cite[Proof of Lemma 2]{NuOL} (which are based
on an application of the multiplication formulae for multiple integrals, see
\cite[Proposition 1.1.3]{NualartBook}), one obtains that%
\begin{eqnarray*}
\left\Vert D\left[ I_{d_{j}}\left( f_{l}^{\left( j\right) }\right) \right]
\right\Vert _{\mathfrak{H}}^{2} &=&d_{j}\left( d_{j}!\right) \left\Vert
f_{l}^{\left( j\right) }\right\Vert _{\mathfrak{H}^{\otimes d_{j}}}^{2} \\
&&+d_{j}^{2}\sum_{p=1}^{d_{j}-1}\left( p-1\right) !\binom{n-1}{p-1}%
^{2}I_{2\left( d_{j}-p\right) }\left( f_{l}^{\left( j\right) }\widetilde{%
\otimes }_{p}f_{l}^{\left( j\right) }\right) \text{.}
\end{eqnarray*}%
Since $\left\Vert f_{l}^{\left( j\right) }\otimes _{p}f_{l}^{\left( j\right)
}\right\Vert _{\mathfrak{H}^{\otimes 2\left( d_{j}-p\right) }}^{2}$ $\geq
\left\Vert f_{l}^{\left( j\right) }\widetilde{\otimes }_{p}f_{l}^{\left(
j\right) }\right\Vert _{\mathfrak{H}^{\otimes 2\left( d_{j}-p\right) }}^{2}$%
, the last relation implies immediately that 1. $\Rightarrow $ 5.. To prove
the opposite implication, first observe that, due to the boundedness of (\ref%
{Avarbound}) and the Cauchy-Schwarz inequality, there exists a finite
constant $M$ (independent of $j$ and $l$) such that%
\begin{equation*}
\left\Vert f_{l}^{\left( j\right) }\otimes _{p}f_{l}^{\left( j\right)
}\right\Vert _{\mathfrak{H}^{\otimes 2\left( d_{j}-p\right) }}^{2}\leq
\left\Vert f_{l}^{\left( j\right) }\right\Vert _{\mathfrak{H}^{\otimes
d_{j}}}^{4}\leq M\text{.}
\end{equation*}%
This implies that, for every sequence $\left\{ l_{n}\right\} $, there exists
a subsequence $\left\{ l_{n^{\prime }}\right\} $ such that the sequences $%
\left\Vert f_{l_{n^{\prime }}}^{\left( j\right) }\otimes
_{p}f_{l_{n^{\prime }}}^{\left( j\right) }\right\Vert
_{\mathfrak{H}^{\otimes 2\left( d_{j}-p\right) }}^{2}$ and
$d_{j}!\left\Vert f_{l_{n^{\prime }}}^{\left( j\right)
}\right\Vert^{2} _{\mathfrak{H}^{\otimes d_{j}}}$ are convergent for every $%
j=1,...,k$ and every $p=1,...,d_{j}-1$ (recall that, by assumption, there
exists a constant $\eta >0$, such that $\left\Vert f_{l_{n^{\prime
}}}^{\left( j\right) }\right\Vert _{\mathfrak{H}^{\otimes d_{j}}}$ $\geq
\eta $, for every $j$ and $l$). We shall now prove that, whenever (\ref{Mall}%
) is verified, then necessarily $\left\Vert f_{l_{n^{\prime }}}^{\left(
j\right) }\otimes _{p}f_{l_{n^{\prime }}}^{\left( j\right) }\right\Vert _{%
\mathfrak{H}^{\otimes 2\left( d_{j}-p\right) }}^{2}\rightarrow 0$.
Indeed, Theorem 4 in \cite{NuOL} implies that, if (\ref{Mall})
takes place and $d_{j}!
\left\Vert f_{l_{n^{\prime }}}^{\left( j\right) }\right\Vert^{2} _{\mathfrak{H}%
^{\otimes d_{j}}}\rightarrow c>0$, then necessarily
\begin{equation}
I_{d_{j}}\left( f_{l_{n^{\prime }}}^{\left( j\right) }\right) \overset{Law}{%
\rightarrow }N\left( 0,c\right) \text{,}  \label{f}
\end{equation}%
where $N\left( 0,c\right) $ stands for a centered Gaussian random variable
with variance $c$. But Theorem 1 in \cite{NuPe} implies that, if (\ref{f})
is verified, then $\left\Vert f_{l_{n^{\prime }}}^{\left( j\right) }\otimes
_{p}f_{l_{n^{\prime }}}^{\left( j\right) }\right\Vert _{\mathfrak{H}%
^{\otimes 2\left( d_{j}-p\right) }}^{2}\rightarrow 0$, thus proving our
claim. This shows that 5. $\Rightarrow $ 1..
\end{proof}

\bigskip

The next result says that, under the additional assumption that the
variances of the elements of $\mathbf{I}\left( l\right) $ converge to one,
the asymptotic approximation (\ref{bb}) is equivalent to the fact that each
component of $\mathbf{I}\left( l\right) $ verifies a CLT. The proof is
elementary, and therefore omitted.

\begin{corollary}
\label{C : NOL}Fix $k\geq 2$, and suppose that the sequence $\mathbf{I}%
\left( l\right) $, $l\geq 1$, is such that, for every $j=1,...,k$, the
sequence of variances appearing in (\ref{Avarbound}) converges to 1, as $%
l\rightarrow +\infty $. Then, each one of Conditions 1.-5. in the statement
of Theorem \ref{T : Gen} is equivalent to the following: for every $%
j=1,...,k,$%
\begin{equation}
I_{d_{j}}\left( f_{l}^{\left( j\right) }\right) \underset{l\rightarrow
+\infty }{\overset{Law}{\rightarrow }}N\left( 0,1\right) \text{,}
\label{cltt}
\end{equation}%
where $N\left( 0,1\right) $ is a centered Gaussian random variable with
unitary variance.
\end{corollary}

\bigskip

\textbf{Remark. }The results of this section can be suitably extended to
deal with the Gaussian approximations of random vectors of the type $%
(F_{l}^{\left( 1\right) }\left( X\right) ,...,F_{l}^{\left( k\right) }\left(
X\right) )$, where $F_{l}^{\left( j\right) }\left( X\right) $, $j=1,...,k$,
is a general square integrable functional of the isonormal process $X$, not
necessarily having the form of a multiple integral. See \cite[Th. 6]{MaPe2}
for a statement containing an extension of this type.

\section{Proof of Theorem \protect\ref{T : Appendix} \label{S : Proof}}

We provide the proof in the case where%
\begin{equation}
\mathfrak{H}=L^{2}\left( \left[ 0,1\right] ,\mathcal{B}\left( \left[ 0,1%
\right] \right) ,dx\right) =L^{2}(\left[ 0,1\right] ),  \label{leb}
\end{equation}
where $dx$ stands for Lebesgue measure. The extension to a general $%
\mathfrak{H}$ is obtained by using the same arguments outlined in \cite[%
Section 2.2]{NuPe}. If $\mathfrak{H}$ is as in (\ref{leb}), then for every $%
d\geq 2$ one has that $\mathfrak{H}^{\odot d}=L_{s}^{2}(\left[ 0,1\right]
^{d})$, where the symbol $L_{s}^{2}(\left[ 0,1\right] ^{d})$ indicates the
class of symmetric, real-valued and square-integrable functions (with
respect to the Lebesgue measure) on $\left[ 0,1\right] ^{d}$. Also, the
isonormal process $X$ coincides with the Gaussian space generated by the
standard Brownian motion
\begin{equation*}
t\mapsto W_{t}\triangleq X\left( 1_{\left[ 0,t\right] }\right) ,\text{ \ \ }%
t\in \left[ 0,1\right] .
\end{equation*}%
This implies in particular that, for every $d\geq 2$, the Wiener-It\^{o}
integral $I_{d}\left( f\right) $, $f\in L_{s}^{2}\left( \left[ 0,1\right]
^{d}\right) $, can be rewritten in terms of an iterated stochastic integral
with respect to $W$, that is:
\begin{equation}
I_{d}\left( f\right) =d!\int_{0}^{1}\int_{0}^{t_{1}}\cdot \cdot \cdot
\int_{0}^{t_{d-1}}f\left( t_{1},...,t_{d}\right) dW_{t_{d}}\cdot \cdot \cdot
dW_{t_{2}}dW_{t_{1}}.  \label{Amwi}
\end{equation}
We also have that $I_{1}\left( f\right) =\int_{0}^{1}f\left(
s\right) dW_{s}$ for every $f\in L_{s}^{2}(\left[ 0,1\right]
^{1})\equiv $ $L^{2}\left( \left[ 0,1\right] \right) $. Note that
the RHS of (\ref{Amwi}) is just an iterated adapted stochastic
integral of the It\^{o} type. Finally, for every $f\in
L_{s}^{2}(\left[ 0,1\right] ^{d})$, every $g\in L_{s}^{2}(\left[
0,1\right] ^{d'})$ and every $p=0,...,d\wedge d'$, we observe that
the contraction $f\otimes _{p}g$ is the (not necessarily
symmetric) element of $L^{2}(\left[ 0,1\right] ^{d+d'-2p})$
given by:%
\begin{eqnarray}
f\otimes _{p}g\left( y_{1},...,y_{d+d'-2p}\right) &=&\int_{ \left[
0,1\right] ^{p}}f\left( y_{1},...,y_{d-p},a_{1},...,a_{p}\right)
\times  \label{Acont} \\
&&\text{ \ \ \ \ \ \ \ \ \ \ \ \ \ \ \ \ \ \ }\times g\left(
y_{d-p+1},...,y_{d+d'-2p},a_{1},...,a_{p}\right) da_{1}...da_{p}%
\text{.}  \notag
\end{eqnarray}%

In the framework of (\ref{leb}), the proof of Theorem \ref{T : Appendix}
relies on some computations contained in \cite{PT}, as well as on an
appropriate use of the \textsl{Burkholder-Davis-Gundy inequalities} (see for
instance \cite[Ch. IV \S 4]{RY}). Fix $\mathbf{\lambda }=\left( \lambda
_{1},...,\lambda _{k}\right) \in \mathbb{R}^{k}$, and consider the random
variable
\begin{eqnarray*}
\left\langle \mathbf{\lambda },\mathbf{I}\left( l\right) \right\rangle _{k}
&=&\sum_{j=1}^{k}\lambda _{j}d_{j}!\int_{0}^{1}\cdot \cdot \cdot
\int_{0}^{u_{d_{j}-1}}f_{l}^{\left( j\right) }\left(
u_{1},...,u_{d_{j}}\right) dW_{u_{d\,j}}\cdot \cdot \cdot dW_{u_{1}} \\
&\triangleq &\sum_{j=1}^{k}\lambda _{j}d_{j}!J_{d_{j}}^{1}\left(
f_{l}^{\left( j\right) }\right) =\int_{0}^{1}\left( \sum_{j=1}^{k}\lambda
_{j}d_{j}!J_{d_{j}-1}^{u}\left( f_{l}^{\left( j\right) }\left( u,\cdot
\right) \right) \right) dW_{u} \\
&=&\int_{0}^{1}\left( \sum_{j=1}^{k}\lambda _{j}d_{j}I_{d_{j}-1}\left(
f_{l}^{\left( j\right) }\left( u,\cdot \right) \mathbf{1}_{\left[ 0,u\right]
^{d_{j}-1}}\right) \right) dW_{u}\text{,}
\end{eqnarray*}%
where, for every $d\geq 1$, every $t\in \left[ 0,1\right] $ and every $f\in
L_{s}^{2}\left( \left[ 0,1\right] ^{d}\right) $, we define $J_{d}^{t}\left(
f\right) =I_{d}\left( f\mathbf{1}_{\left[ 0,t\right] ^{d}}\right) /d!$ (for
every $c\in \mathbb{R}$, we also use the conventional notation $%
J_{0}^{t}\left( c\right) =c$). We start by recalling some preliminary
results involving Brownian martingales. Start by setting, for every $u\in %
\left[ 0,1\right] $, $\phi _{\mathbf{\lambda },l}\left( u\right)
=\sum_{j=1}^{k}\lambda _{j}d_{j}I_{d_{j}-1}\left( f_{l}^{\left( j\right)
}\left( u,\cdot \right) \mathbf{1}_{\left[ 0,u\right] ^{d_{j}-1}}\right) $,
and observe that the random application
\begin{equation*}
t\mapsto \sum_{j=1}^{k}\lambda _{j}d_{j}!J_{d_{j}}^{t}\left( f_{l}^{\left(
j\right) }\right) =\int_{0}^{t}\phi _{\mathbf{\lambda },l}\left( u\right)
dW_{u}\text{, \ \ }t\in \left[ 0,1\right] \text{,}
\end{equation*}%
defines a (continuous) square-integrable martingale started from zero, with
respect to the canonical filtration of $W$, noted $\left\{ \mathcal{F}%
_{t}^{W}:t\in \left[ 0,1\right] \right\} $. The quadratic variation of this
martingale is classically given by $t\mapsto \int_{0}^{t}\phi _{\mathbf{%
\lambda },l}\left( u\right) ^{2}du$, and a standard application of the
Dambis, Dubins and Schwarz Theorem (see \cite[Ch. V \S 1]{RY}) yields that,
for every $l\geq 1$, there exists a standard Brownian motion (initialized at
zero) $W^{\left( \mathbf{\lambda },l\right) }=\left\{ W_{t}^{\left( \mathbf{%
\lambda },l\right) }:t\geq 0\right\} $ such that
\begin{equation*}
\left\langle \mathbf{\lambda },\mathbf{I}\left( l\right) \right\rangle
_{k}=\int_{0}^{1}\phi _{\mathbf{\lambda },l}\left( u\right)
dW_{u}=W_{\int_{0}^{1}\phi _{\mathbf{\lambda },l}\left( u\right)
^{2}du}^{\left( \mathbf{\lambda },l\right) }\text{.}
\end{equation*}%
Note that, in general, the definition of $W^{\left( \mathbf{\lambda }%
,l\right) }$ strongly depends on $\mathbf{\lambda }$ and $l$, and that $%
W^{\left( \mathbf{\lambda },l\right) }$ \textsl{is not} a $\mathcal{F}%
_{t}^{W}$-Brownian motion. However, the following relation links the two
Brownian motions $W^{\left( \mathbf{\lambda },l\right) }$ and $W$: there
exists a (continuous) filtration $\left\{ \mathcal{G}_{t}^{\left( \mathbf{%
\lambda },l\right) }:t\geq 0\right\} $ such that (i) $W_{t}^{\left( \mathbf{%
\lambda },l\right) }$ is a $\mathcal{G}_{t}^{\left( \mathbf{\lambda }%
,l\right) }$-Brownian motion, and (ii) for every fixed $s\in \left[ 0,1%
\right] $ the positive random variable $\int_{0}^{s}\phi _{\mathbf{\lambda }%
,l}\left( u\right) ^{2}du$ is a $\mathcal{G}_{t}^{\left( \mathbf{\lambda }%
,l\right) }$-stopping time. Now define the positive constant (which is
trivially a $\mathcal{G}_{t}^{\left( \mathbf{\lambda },l\right) }$-stopping
time)
\begin{equation*}
q\left( \mathbf{\lambda },l\right) =\int_{0}^{1}\mathbb{E(}\phi _{\mathbf{%
\lambda },l}\left( u\right) ^{2})du\text{,}
\end{equation*}%
and observe that the usual properties of complex exponentials and a standard
application of the Burkholder-Davis-Gundy inequality (in the version stated
in \cite[Corollary 4.2, Ch. IV ]{RY}) yield the following estimates:%
\begin{eqnarray}
\left\vert \mathbb{E}\left[ \exp \left( \mathrm{i}\left\langle \mathbf{%
\lambda },\mathbf{I}\left( l\right) \right\rangle _{k}\right) \right] -%
\mathbb{E}\left[ \exp \left( \mathrm{i}W_{q\left( \mathbf{\lambda },l\right)
}^{\left( \mathbf{\lambda },l\right) }\right) \right] \right\vert
&=&\left\vert \mathbb{E}\left[ \exp \left( \mathrm{i}W_{\int_{0}^{1}\phi _{%
\mathbf{\lambda },l}\left( u\right) ^{2}du}^{\left( \mathbf{\lambda }%
,l\right) }\right) \right] -\mathbb{E}\left[ \exp \left( \mathrm{i}%
W_{q\left( \mathbf{\lambda },l\right) }^{\left( \mathbf{\lambda },l\right)
}\right) \right] \right\vert   \notag \\
&\leq &\mathbb{E}\left[ \left\vert W_{\int_{0}^{1}\phi _{\mathbf{\lambda }%
,l}\left( u\right) ^{2}du}^{\left( \mathbf{\lambda },l\right) }-W_{q\left(
\mathbf{\lambda },l\right) }^{\left( \mathbf{\lambda },l\right) }\right\vert %
\right]   \notag \\
&\leq &\mathbb{E}\left[ \left\vert W_{\int_{0}^{1}\phi _{\mathbf{\lambda }%
,l}\left( u\right) ^{2}du}^{\left( \mathbf{\lambda },l\right) }-W_{q\left(
\mathbf{\lambda },l\right) }^{\left( \mathbf{\lambda },l\right) }\right\vert
^{4}\right] ^{\frac{1}{4}}  \notag \\
&\leq &C\mathbb{E}\left[ \left\vert \int_{0}^{1}\phi _{\mathbf{\lambda }%
,l}\left( u\right) ^{2}du-q\left( \mathbf{\lambda },l\right) \right\vert ^{2}%
\right] ^{\frac{1}{4}}\text{,}  \label{Aest}
\end{eqnarray}%
where $C$ is some universal constant independent of $\mathbf{\lambda }$ and $%
l$. To see how to obtain the inequality (\ref{Aest}), introduce first the
shorthand notation $T\left( \mathbf{\lambda },l\right) \triangleq
\int_{0}^{1}\phi _{\mathbf{\lambda },l}\left( u\right) ^{2}du$ (recall that $%
T\left( \mathbf{\lambda },l\right) $ is a $\mathcal{G}_{t}^{\left( \mathbf{%
\lambda },l\right) }$-stopping time), and then write%
\begin{equation*}
\left\vert W_{\int_{0}^{1}\phi _{\mathbf{\lambda },l}\left( u\right)
^{2}du}^{\left( \mathbf{\lambda },l\right) }-W_{q\left( \mathbf{\lambda }%
,l\right) }^{\left( \mathbf{\lambda },l\right) }\right\vert =\left\vert
\int_{T\left( \mathbf{\lambda },l\right) \wedge q\left( \mathbf{\lambda }%
,l\right) }^{T\left( \mathbf{\lambda },l\right) \vee q\left( \mathbf{\lambda
},l\right) }dW_{u}^{\left( \mathbf{\lambda },l\right) }\right\vert
=\left\vert \int_{0}^{T\left( \mathbf{\lambda },l\right) \vee q\left(
\mathbf{\lambda },l\right) }H\left( u\right) dW_{u}^{\left( \mathbf{\lambda }%
,l\right) }\right\vert \text{,}
\end{equation*}%
where $H\left( u\right) $ is the $\mathcal{G}_{u}^{\left( \mathbf{\lambda }%
,l\right) }$-predictable process given by $H\left( u\right) =\mathbf{1}%
\left\{ u\geq T\left( \mathbf{\lambda },l\right) \wedge q\left( \mathbf{%
\lambda },l\right) \right\} $, so that
\begin{eqnarray*}
\left\vert \int_{0}^{T\left( \mathbf{\lambda },l\right) \vee q\left( \mathbf{%
\lambda },l\right) }H\left( u\right) ^{2}du\right\vert  &=&\left\vert
T\left( \mathbf{\lambda },l\right) \wedge q\left( \mathbf{\lambda },l\right)
-T\left( \mathbf{\lambda },l\right) \vee q\left( \mathbf{\lambda },l\right)
\right\vert  \\
&=&\left\vert T\left( \mathbf{\lambda },l\right) -q\left( \mathbf{\lambda }%
,l\right) \right\vert =\left\vert \int_{0}^{1}\phi _{\mathbf{\lambda }%
,l}\left( u\right) ^{2}du-q\left( \mathbf{\lambda },l\right) \right\vert
\text{.}
\end{eqnarray*}%
In particular, relation (\ref{Aest}) yields that the proof of Theorem \ref{T
: Appendix} is concluded, once the following two facts are proved: (A) $%
W_{q\left( \mathbf{\lambda },l\right) }^{\left( \mathbf{\lambda },l\right)
}=\left\langle \mathbf{\lambda },\mathbf{N}\left( l\right) \right\rangle _{k}
$, for every $\mathbf{\lambda }\in \mathbb{R}^{k}$ and every $l\geq 1$; (B)
the sequence
\begin{equation*}
\mathbb{E}\left[ \left\vert \int_{0}^{1}\phi _{\mathbf{\lambda },l}\left(
u\right) ^{2}du-q\left( \mathbf{\lambda },l\right) \right\vert ^{2}\right]
\text{, \ \ }l\geq 1\text{,}
\end{equation*}%
converges to zero, uniformly in $\mathbf{\lambda }$, on every compact set of
the type $M=\left[ -T,T\right] ^{k}$, where $T\in \left( 0,+\infty \right) $%
. The proof of (A) is immediate: indeed, $W^{\left( \mathbf{\lambda }%
,l\right) }$ is a standard Brownian motion and, by using the isometric
properties of stochastic integrals and the fact that the covariance
structures of $\mathbf{N}\left( l\right) $ and $\mathbf{I}\left( l\right) $
coincide,
\begin{equation*}
q\left( \mathbf{\lambda },l\right) =\int_{0}^{1}\mathbb{E(}\phi _{\mathbf{%
\lambda },l}\left( u\right) ^{2})du=\mathbb{E}\left[ \left( \int_{0}^{1}\phi
_{\mathbf{\lambda },l}\left( u\right) dW_{u}\right) ^{2}\right] =\mathbb{E}%
\left[ \left\langle \mathbf{\lambda },\mathbf{I}\left( l\right)
\right\rangle _{k}^{2}\right] =\mathbb{E}\left[ \left\langle \mathbf{\lambda
},\mathbf{N}\left( l\right) \right\rangle _{k}^{2}\right] .
\end{equation*}%
To prove (B), use a standard version of the multiplication formula between
multiple stochastic integrals (see for instance \cite[Proposition 1.5.1]%
{NualartBook})
\begin{eqnarray}
&&\int_{0}^{1}\phi _{\mathbf{\lambda },l}\left( u\right)
^{2}du=\int_{0}^{1}\left( \sum_{j=1}^{k}\lambda _{j}d_{j}I_{d_{j}-1}\left(
f_{l}^{\left( j\right) }\left( u,\cdot \right) \mathbf{1}_{\left[ 0,u\right]
^{d_{j}-1}}\right) \right) ^{2}du  \notag \\
&=&\int_{0}^{1}\sum_{j,i=1}^{k}\lambda _{j}\lambda
_{i}d_{j}d_{i}I_{d_{i}-1}\left( f_{l}^{\left( i\right) }\left( u,\cdot
\right) \mathbf{1}_{\left[ 0,u\right] ^{d_{i}-1}}\right) I_{d_{j}-1}\left(
f_{l}^{\left( j\right) }\left( u,\cdot \right) \mathbf{1}_{\left[ 0,u\right]
^{d_{j}-1}}\right) du  \notag \\
&=&q\left( \mathbf{\lambda },l\right) +\sum_{j,i=1}^{k}\lambda _{j}\lambda
_{i}d_{j}d_{i}\int_{0}^{1}\sum_{p=0}^{D\left( i,j\right) }\binom{d_{i}-1}{p}%
\binom{d_{j}-1}{p}  \label{Avardev} \\
&&\text{ \ \ \ \ \ \ \ \ \ \ \ \ \ \ \ }\times I_{d_{i}+d_{j}-2-2p}\left(
(f_{l}^{\left( i\right) }\left( u,\cdot \right) \mathbf{1}_{\left[ 0,u\right]
^{d_{i}-1}})\otimes _{p}(f_{l}^{\left( j\right) }\left( u,\cdot \right)
\mathbf{1}_{\left[ 0,u\right] ^{d_{j}-1}})\right) \text{,}  \notag
\end{eqnarray}%
where the index $D\left( i,j\right) $ is defined as%
\begin{equation*}
D\left( i,j\right) =\left\{
\begin{array}{ll}
d_{i}-2 & \text{if }d_{i}=d_{j} \\
\min \left( d_{i},d_{j}\right) -1 & \text{if \ }d_{i}\neq d_{j}\text{.}%
\end{array}%
\right.
\end{equation*}%
Formula (\ref{Avardev}) implies that, for every $\mathbf{\lambda }\in \left[
-T,T\right] ^{k}$ ($T>0$)$,$%
\begin{eqnarray}
&&\mathbb{E}\left[ \left\vert \int_{0}^{1}\phi _{\mathbf{\lambda },l}\left(
u\right) ^{2}du-q\left( \mathbf{\lambda },l\right) \right\vert ^{2}\right] ^{%
\frac{1}{2}}  \notag \\
&\leq &(T\max_{i}d_{i})^{2}\sum_{i,j=1}^{k}\sum_{p=0}^{D\left( i,j\right) }%
\binom{d_{i}-1}{p}\binom{d_{j}-1}{p}  \label{Arhs} \\
&&\times \mathbb{E}\left[ \left( \int_{0}^{1}I_{d_{i}+d_{j}-2-2p}\left(
(f_{l}^{\left( i\right) }\left( u,\cdot \right) \mathbf{1}_{\left[ 0,u\right]
^{d_{i}-1}})\otimes _{p}(f_{l}^{\left( j\right) }\left( u,\cdot \right)
\mathbf{1}_{\left[ 0,u\right] ^{d_{j}-1}})\right) du\right) ^{2}\right] ^{%
\frac{1}{2}}  \notag
\end{eqnarray}%
(note that the RHS of (\ref{Arhs}) does not depend on $\mathbf{\lambda }$).
Finally, a direct application of the calculations contained in \cite[p.
253-255]{PT} yields that, for every $i,j=1,...,k$ and every $p=0,...,D\left(
i,j\right) $,%
\begin{equation}
\mathbb{E}\left[ \left( \int_{0}^{1}I_{d_{i}+d_{j}-2-2p}\left(
(f_{l}^{\left( i\right) }\left( u,\cdot \right) \mathbf{1}_{\left[ 0,u\right]
^{d_{i}-1}})\otimes _{p}(f_{l}^{\left( j\right) }\left( u,\cdot \right)
\mathbf{1}_{\left[ 0,u\right] ^{d_{j}-1}})\right) du\right) ^{2}\right] ^{%
\frac{1}{2}}\rightarrow 0\text{,}  \label{q}
\end{equation}%
as $l\rightarrow +\infty $. This concludes the proof of Theorem \ref{T :
Appendix}. \ $\blacksquare $

\bigskip

\textbf{Remark. }By inspection of the calculations contained in \cite[p.
253-255]{PT}, it is easily seen that, to deduce (\ref{q}) from (\ref%
{Acontrconv}), it is necessary that the sequence of variances (\ref%
{Avarbound}) is bounded.

\section{Concluding remarks on applications\label{S : Concl}}

Theorem \ref{T : Appendix} and Theorem \ref{T : Gen} are used in \cite{MaPe2}
to deduce high-frequency asymptotic results for subordinated spherical
random fields. This study is strongly motivated by the probabilistic
modelling and statistical analysis of the Cosmic Microwave Background
radiation (see \cite{Marinucci}, \cite{MarPTRF}, \cite{MaPe} and \cite{MaPe2}
for a detailed discussion of these applications). In what follows, we
provide a brief presentation of some of the results obtained in \cite{MaPe2}.

Let $\mathbb{S}^{2}=\left\{ x\in \mathbb{R}^{3}:\left\Vert x\right\Vert
=1\right\} $ be the unit sphere, and let $T=\{T\left( x\right) :x\in \mathbb{%
S}^{2}\}$ be a real-valued (centered) Gaussian field which is also \textsl{%
isotropic}, in the sense that $T\left( x\right) \overset{Law}{=}T\left(
\mathcal{R}x\right) $ (in the sense of stochastic processes) for every
rotation $\mathcal{R}\in SO\left( 3\right) $. The following facts are well
known:

\begin{description}
\item[(1)] The trajectories of $T$ admit the harmonic expansion $T\left(
x\right) $ $=\sum_{l=0}^{\infty }\sum_{m=-l}^{l}a_{lm}Y_{lm}\left( x\right) $%
, where $\{Y_{lm}:l\geq 0$, \ $m=-l,...,l\}$ is the class of \textsl{%
spherical harmonics} (defined e.g. in \cite[Ch. 5]{VMK});

\item[(2)] The complex-valued array of harmonic coefficients $\{a_{lm}:l\geq
0$, $l\geq 0,$\ $m=-l,...,l\}$ is composed of centered Gaussian random
variables such that the variances $\mathbb{E}\left\vert a_{lm}\right\vert
^{2}$ $\triangleq C_{l}$ depend exclusively on $l$ (see for instance \cite%
{BaMa});

\item[(3)] The law of $T$ is completely determined by the \textsl{power
spectrum }$\left\{ C_{l}:l\geq 0\right\} $ defined at the previous point.
\end{description}

Now fix $q\geq 2$, and consider the subordinated field
\begin{equation*}
T^{\left( q\right) }\left( x\right) \triangleq H_{q}\left( T\left( x\right)
\right) \text{, \ \ }x\in \mathbb{S}^{2}\text{,}
\end{equation*}%
where $H_{q}$ is the $q$th Hermite polynomial. Plainly, the field $T^{\left(
q\right) }$ is isotropic and admits the harmonic expansion%
\begin{equation*}
T^{\left( q\right) }\left( x\right) =\sum_{l=0}^{\infty
}\sum_{m=-l}^{l}a_{lm;q}Y_{lm}\left( x\right) \triangleq \sum_{l=0}^{\infty
}T_{l}^{\left( q\right) }\left( x\right) \text{,}
\end{equation*}%
where $a_{lm;q}$ $\triangleq $ $\int_{\mathbb{S}^{2}}T^{\left( q\right)
}\left( z\right) \overline{Y_{lm}\left( z\right) }dz$. For every $l\geq 0$,
the field $T_{l}^{\left( q\right) }=\sum_{m=-l}^{l}a_{lm;q}Y_{lm}$ is
real-valued and isotropic, and it is called the $l$th \textsl{frequency
component} of $T^{\left( q\right) }$ (see \cite{Marinucci} or \cite{MaPe2}
for a physical interpretation of frequency components). In \cite{MaPe2}, the
following problem is studied.

\bigskip

\textbf{Problem A. }\textsl{Fix }$q\geq 2$\textsl{. Find conditions on the
power spectrum }$\left\{ C_{l}:l\geq 0\right\} $\textsl{\ to have that the
finite dimensional distributions (f.d.d.'s) of the normalized frequency field%
}%
\begin{equation*}
\overline{T}_{l}^{\left( q\right) }\left( x\right) \triangleq \frac{%
T_{l}^{\left( q\right) }\left( x\right) }{\mathbf{Var}\left( T_{l}^{\left(
q\right) }\left( x\right) \right) ^{1/2}}\text{, \ \ }x\in \mathbb{S}^{2}%
\text{,}
\end{equation*}%
\textsl{are `asymptotically close to Gaussian', as }$l\rightarrow +\infty $.

\bigskip

The main difficulty when dealing with Problem A is that (due to isotropy)
one has always that%
\begin{equation}
\mathbb{E}\left[ \overline{T}_{l}^{\left( q\right) }\left( x\right)
\overline{T}_{l}^{\left( q\right) }\left( y\right) \right] =P_{l}\left( \cos
\left\langle x,y\right\rangle \right) \text{,}  \label{cobhat}
\end{equation}%
where $P_{l}$ is the $l$th Legendre polynomial, and $\left\langle
x,y\right\rangle $ is the angle between $x$ and $y$. Indeed, since in
general the quantity $P_{l}\left( \cos \left\langle x,y\right\rangle \right)
$ does not converge (as $l\rightarrow +\infty $), one cannot prove that the
f.d.d.'s of $\overline{T}_{l}^{\left( q\right) }$ converge to those of a
Gaussian field (even if $\overline{T}_{l}^{\left( q\right) }\left( x\right) $
converges in law to a Gaussian random variable for every fixed $x$).
However, as an application of Theorem \ref{T : Appendix} and Proposition \ref%
{C : Appendix}, one can prove the following approximation result.

\begin{proposition}
Under the above notation and assumptions, suppose that, for any fixed $x\in
\mathbb{S}^{2}$,
\begin{equation}
\overline{T}_{l}^{\left( q\right) }\left( x\right) \underset{l\rightarrow
+\infty }{\overset{Law}{\rightarrow }}N\left( 0,1\right) \text{.}
\label{ccff}
\end{equation}%
Then, for any $k\geq 1$, any $x_{1},...,x_{k}\in \mathbb{S}^{2}$ and any $%
\gamma $ metrizing the weak convergence on $\mathbf{P}\left( \mathbb{R}%
^{k}\right) $,
\begin{equation}
\gamma \left( \mathcal{L}\left( \overline{T}_{l}^{\left( q\right) }\left(
x_{1}\right) ,...,\overline{T}_{l}^{\left( q\right) }\left( x_{k}\right)
\right) ,\mathbf{N}\left( l\right) \right) \underset{l\rightarrow +\infty }{%
\rightarrow }0\text{,}  \label{hol}
\end{equation}%
where, for every $l$, $\mathbf{N}\left( l\right) =\left(
N_{l}^{\left( 1\right) },...,N_{l}^{\left( k\right) }\right) $ is
a centered real-valued Gaussian vector such that
\begin{equation*}
\mathbb{E}\left\{ N_{l}^{\left( i\right) }N_{l}^{\left( j\right)
}\right\} =P_{l}\left( \cos \left\langle x_{i},x_{j}\right\rangle
\right) .
\end{equation*}
\end{proposition}

\begin{proof}
Since $\overline{T}_{l}^{\left( q\right) }\left( x\right) $ is a
linear functional involving uniquely Hermite polynomials of order
$q$ (written on the Gaussian field $T$) one deduces that there
exists a real Hilbert space $\mathfrak{H}$
such that (in the sense of stochastic processes)%
\begin{equation*}
\overline{T}_{l}^{\left( q\right) }\left( x\right) \overset{Law}{=}%
I_{q}\left( f_{\left( q,l,x\right) }\right) \text{,}
\end{equation*}%
where the class of symmetric kernels
\begin{equation*}
\left\{ f_{\left( q,l,x\right) }:l\geq 0,\text{ }x\in \mathbb{S}^{2}\right\}
\end{equation*}%
is a subset of $\mathfrak{H}^{\odot q}$ , and $I_{q}\left( f_{\left(
q,l,x\right) }\right) \ $stands for the $q$th Wiener-It\^{o} integral of $%
f_{\left( q,l,x\right) }$ with respect to an isonormal Gaussian process over
$\mathfrak{H}$, as defined in Section \ref{S : Pre}. Since the variances of
the components of the vector $(\overline{T}_{l}^{\left( q\right) }\left(
x_{1}\right) ,...,\overline{T}_{l}^{\left( q\right) }\left( x_{k}\right) )$
are all equal to 1 by construction, we can apply Theorem \ref{T : Gen} and
Proposition \ref{C : Appendix}. Indeed, by Theorem \ref{T : Gen} we know
that (\ref{ccff}) implies that, for every $p=1,...,q-1$ and every $%
j=1,...,k, $%
\begin{equation*}
f_{\left( q,l,x_{j}\right) }\otimes _{p}f_{\left( q,l,x_{j}\right)
}\rightarrow 0\text{ in }\mathfrak{H}^{\odot 2\left( q-p\right) }.
\end{equation*}%
Finally, Proposition \ref{C : Appendix} and (\ref{cobhat}) imply immediately
the desired conclusion.
\end{proof}

\bigskip

The derivation of sufficient conditions to have (\ref{ccff}) is the main
object of \cite{MaPe2}. In particular, it is proved that sufficient (and
sometimes also necessary) conditions for (\ref{ccff}) can be neatly
expressed in terms of the so-called \textsl{Clebsch-Gordan} \textsl{%
coefficients }(see again \cite{VMK}), that are elements of unitary matrices
connecting reducible representations of $SO\left( 3\right) $.

\bigskip

\textbf{Acknowledgements -- }I\ am grateful to D.\ Marinucci for many
fundamental discussions on the subject of this paper. Part of this work has
been written when I was visiting the Departement of Statistics and Applied
Mathematics of Turin University. I\ wish to thank M. Marinacci and I. Pr\"{u}%
nster for their hospitality.

\end{document}